\documentclass[12pt]{article}
\usepackage{amssymb, amsmath, amsthm, bbm, verbatim}
\newcommand{\ravno}{\mathop{=}}
\newcommand{\too}{\mathop{\longrightarrow}}
\newcommand{\supp}{\mbox{supp}\,}

\newtheorem{teo}{Theorem}
\newtheorem{prop}{Proposition}

\newtheorem{cor}{Corollary}
\newtheorem{rem}{Remark}
\newtheorem{lem}{Lemma}

\title{Walsh and wavelet methods for\\ differential equations on the Cantor group
}
\author{
E.  Lebedeva
\footnote{Mathematics and Mechanics Faculty, St. Petersburg State University,
Universitetsky prospekt, 28, Peterhof,  Saint Petersburg,
 198504, Russia; Institute of Applied Mathematics and Mechanics, St. Petersburg State Polytechnical University,
 Polytechnicheskay 29, 195251, St. Petersburg, Russia},
 M. Skopina
\footnote{Faculty of Applied Mathematics and Control Processes, St. Petersburg State University,
Universitetsky prospekt, 35, Peterhof, Saint Petersburg,   198504,
Russia}
}
\date{
 ealebedeva2004@gmail.com, skopina@ms1167.spb.edu}

\begin{document}
\maketitle

\newcommand{\nul}{{\bf0}}
\newcommand{\rd}{{\mathbb R}^d}
\newcommand{\zd}{{\mathbb Z}^{d}}
\renewcommand{\r}{{\mathbb R}}
\newcommand{\z} {{\mathbb Z}}
\newcommand{\cn} {{\mathbb C}}
\newcommand{\n} {{\mathbb N}}

\begin{abstract}
Ordinary and partial differential equation for unknown functions defined
on the Cantor dyadic  group are studied. We consider two types of equations:
related to the Gibbs derivatives and to the fractional modified Gibbs
derivatives (or pseudo differential-operators). We find solutions in  classes of
distributions and study under what assumptions these solutions are regular functions with some ''good'' properties.
\end{abstract}



\section{Introduction}

Dyadic (Walsh) analysis is actively studied for the last three decades.
The foundations of the Walsh theory are presented in monographs~\cite{GES}
and~\cite{SWS}. The theory may be stated in two equivalent forms:
the functions under consideration are defined on the real semiline or
they are defined on the Cantor group (see~\cite[\S 1.2]{GES}).
We work  in the latter form.

The dyadic wavelet theory is actively studied nowadays.
In 1996 the concept  of   multiresolution analysis (MRA) for the Cantor  group
was introduced by Lang~\cite{Lang}, who also developed
 a method for the construction MRA-based orthogonal wavelet bases.
Later  deep investigations of the foundation of the dyadic wavelet theory
were added by Farkov and Protasov~\cite{FP}.


A concept of dyadic derivative was introduced by Gibbs~\cite{Gibbs} in 1967.
Numerous generalizations of this notion in different directions can be found
in the literature, the surveys on this topic with extensive bibliographies were given in~\cite{SG}, \cite{SA}. A notion of so-called modified Gibbs derivative we use in the present paper was introduced by Golubov in~\cite{G}.

Dyadic analogues of the classical partial differential equations were considered by Butzer and Wagner \cite{BW1}. They found a dyadic analogue of d'Alamber's solution to  the  one-dimensional homogeneous wave equation $\Phi^{[2]}_{x^2}=\Phi^{[2]}_{t^2}$,
where $f^{[2]}$ denotes the Gibbs derivative of the second order.
Note that both variables  $x$ and  $t$ are elements of the group $G$ here. In the present paper, we consider PDE's, where the variable $x$ is in $G$
and the variable $t$, interpreted as time, is real.

We study ordinary and partial differential equations with respect to
both Gibbs and modified Gibbs derivatives.
In the case of Gibbs derivative,
we consider some linear differential equations and the Cauchy problem
for the homogeneous wave equation. The Walsh functions are
eigenfunctions for the Gibbs differentiation operator. Due to this property,
using Walsh expansions, we find all solutions of the equations in the class of periodic distributions.  Also we investigate under what conditions the solution is a regular function with "good" properties (belonging to $L_2$ or continuity).

Unfortunately, the Gibbs derivative has the following drawback.
Consider a simple equation, say $f^{[1]}=g$, where $g$ is $1$-periodic.
To solve the solution one  finds the Fourier-Walsh coefficients of $f^{[1]}$,
 which allows to restore~$f$. However if $g$ is restricted to the fundamental
 domain $I$, then one can repeat this trick  using the
 Fourier-Walsh transform instead of the Fourier-Walsh
 coefficients,  and then restore $f$. It appears that such a function $f$ is not a compactly supported and it does not coincide  on $I$ with  the periodic solution.
 The modified Gibbs differentiability is "more local"
property of functions, an analog of the described  drawback  does not hold for the corresponding differential equations. Namely, if $g$ is a $1$-periodic function, $f$ is an unknown function,  then the periodic solution of the equation
${\cal D}f = g$  coincides on $I$ with the solution of the equation
${\cal D}f = g \mathbbm{1}_I,$ where $\mathbbm{1}_I$ is the characteristic function of $I.$
Moreover, the solution of the later equation ${\cal D}f = g \mathbbm{1}_I$ has the same support as the right-hand side.

We consider ordinary and partial differential equations with respect to the fractional modified Gibbs derivative ${\cal D}^{\alpha}$, $\alpha\in\r$.
In particular, we consider the Cauchy problem  for the one-dimensional non-homogeneous heat equation.
It appeared that all elements of the Haar basis are eigenfunctions for
${\cal D}^{\alpha}$. This allows to solve equations using Haar expansions.
We find all solutions on a class of distributions and investigate
under what conditions the solution is a regular function with "good"
properties. The multiresolution structure of the Haar basis is essentially used
in our technique.

The paper is organized as follows. First, we introduce  necessary notations
and basic facts of the Walsh analysis on the Cantor group.
In Subsection~\ref{Haar} we discuss the Haar basis, its
multiresolution structure, Haar and quasi-Haar expansions.
In Section~\ref{classes} we introduce classes of distributions on the
Cantor group and establish some connections between them.
It Section~\ref{Gde} we study differential equation
with respect to the Gibbs derivative.
It Section~\ref{MGd} we study differential equation
with respect to the fractional modified Gibbs derivative.

\section{Notations and Auxiliary Results}
\label{Not}
\subsection{The Cantor dyadic group}
\label{cant_gr}
Here and in what follows, we  use  basic facts
on Cantor group and Walsh analysis from~\cite{GES} and~\cite{SWS}.

As usual,  by $\n$, $\z$, $\z_+$, $\r$ and $\cn$  we denote
the sets of positive integers, integers, non-negative integers,
real numbers and complex numbers respectively.

Let $G$ be the set of  sequences
$x=(x_k)_{k\in\mathbb{Z}},$
where
$x_k\in\{0,\,1\}$
and either there exists $N(x)\in\mathbb{Z}$
such that $x_{N(x)}=1 $ and $x_k=0$ for $k<N(x)$ or
$x_k=0$ for  all $k\in\mathbb{Z}$ (the latter element $x$ we denote by ${\bf 0}$).
We say that $x$  terminates with  1's (0's) if $x_k=0$ ($x_k=1$) only for finitely many $k$. Supply $G$ with coordinate-wise addition, i.e. define
{the sum} of $x\in G$  and $y\in G$   by
$$
x \oplus y :=(|x_k-y_k|)_{k\in \mathbb{Z}}.
$$
Then $(G,\,\oplus)$  is an abelian group called \texttt{the Cantor dyadic group}.

To supply this group with a metric $d$ we introduce
a map $\|\cdot\|:\ G\to [0,+\infty)$
defined by $\|{\bf 0}\|:=0$ and $\|x\|:=2^{-N(x)}$  for $x\ne{\bf 0}$,
and let $d(x,y):=\|x\oplus y\|$ for  $x,y\in G$. Evidently, $d$ is
a  metric. Moreover, the strong triangle inequality
$d(x,y)\leq \max \{d(x,z),\,d(y,z)\}$ is satisfied, so
$d$ is a non-Archimedean metric. Given $n\in\z$ and $x\in G$,
denote by $I_{n}(x)$
the ball of radius $2^{-n}$ with the center at  $x$, i.e.
$$
I_{n}(x)=\{y\in G: d(x,y) < 2^{-n}\}.
$$
For convenience we set $I_j:=I_j({\bf 0})$ and $I:=I_0$. Note that $I$ is a subgroup of the group $G$.
Since the metric is non-Archimedean, any two balls in $G$ either are disjoint or one
contains the other, and any point of a ball
is the center of the ball. The topological group $G$,
with the topology generated by $d$,
is locally compact and totally disconnected.

Let $\lambda(x):=\sum_{j \in \mathbb{Z}}x_j 2^{-j-1},$ then the map $x\mapsto \lambda(x)$ is a one-to-one correspondence taking $G\setminus\mathbb{Q}_0$ onto $[0,\,\infty),$ where
$\mathbb{Q}_0$ consists of all elements terminating with $1$'s.

We now define dilation on $G$ by $D:\ G\to G$,  where $(Dx)_k=x_{k+1}$
for $x\in G$. We let $D^{-1}:\ G\to G$
 be the inverse mapping  $(D^{-1}x)_k=x_{k-1}$.  We let $D^k$ be
 $D\circ\dots\circ D$
($k$ times) if $k > 0$, and $D^{-1}\circ\dots\circ D^{-1}$ ($-k$ times) if $k < 0$;
$D^0$ is the identity mapping.

We consider functions taking $G$ to $\cn$. Let
$\mathbbm{1}_E$ denote the characteristic function of a set  $E\subset G$.
Note that if $E$ is a ball, then $\mathbbm{1}_E$ is continuous.
For a function
$f:\ G\to \cn$ and a number $h\ge0$, define
$f_{0, h}:\ G\to \cn$ as follows:   for every $x\in G$ we set
$f_{0, h}(x)=f(x\oplus \lambda^{-1}(h))$.
If additionally $j\in\z$, then we set
$$
f_{j,h}(x)=2^{j/2}f_{0,h}(D^jx), \quad x\in G.
$$
We say that a function $f$ defined on $G$ is 1-periodic if
$f(x)=f_{0,1}(x)$ for any $x\in G$. Denote by $C^P$ the space
of 1-periodic continuous functions.

Since $G$ is a locally compact group, there exists the Haar measure
$dx$ on it (see~\cite{HR}), which is  positive,
invariant under the shifts, i.e., $d(x \oplus a)=dx$, and normalized by
$\int_{G}\mathbbm{1}_I(x)\,dx=1$. Hence the functional spaces
$L_q(G)$ and $L_q(E)$, where $E$ is a measurable subset of $G$,
are defined. Denote by $L_q^P$ the space of 1-periodic functions
whose restriction to $I$ is in $L_q(I)$.

Given $t\in G$, the function
$$
\chi_t(x):=(-1)^{\sum_{j \in \mathbb{Z}}t_k x_{-k-1}},\quad x\in G.
$$
is  a group character on $G$. The
 Pontryagin dual group $G^*$ of $G$ is topologically isomorphic to $G$,
 where the isomorphism is given by $t\to\chi_t$.
In the sequel we identify these groups and write $G$ instead of $G^*$.

\texttt{The Fourier-Walsh transform }of $f\in L_1(G)$ is defined by
$$
\widehat{f}(\xi):= \int\limits_{G} f(x) {\rm w}(\xi,\,x)\,dx,
$$
 where  ${\rm w}(t,\,x):=\chi_t(x)$,   $t, x\in G$.
 The Fourier transform is extended to ${ L}_2(G)$ in a
standard way, and the Plancherel equality holds
$$
\langle f,g\rangle:=
\int\limits_{G}f(x)\overline{g(x)}\,dx=
\int\limits_{G}\widehat f(\xi)\overline{\widehat g(\xi)}\,d\xi
=\langle \widehat f, \widehat g\rangle,\quad f,g\in L_2(G).
$$
The inversion formula
\begin{equation}
\label{12}
\widehat{\widehat f}=f
\end{equation}
holds true for any $f\in L_2(G)$.

 Given $n\in\z_+$, set ${\rm w}_n(x):={\rm w}(\lambda^{-1}(n),x)$.
The functions ${\rm w}_n$ are called \texttt{the Walsh functions}.
 These functions are  continuous on $G$ and 1-periodic,
 also they form an orthonormal basis for  $L_2^P$. Thus any $f\in L_2^P$
can be decomposed as
$$
f\ravno^{L_2}\sum_{k=0}^\infty c_{k}{\rm w}_k, \quad c_k\in\cn,
$$
(\texttt{Walsh representation }in the sequel). Here $c_k$ is the $k$-th
Walsh-Fourier coefficient of $f$:
$$
c_k=\widehat f(k):=\int\limits_I f(x){\rm w}_k(x)\,dx.
$$
 A function $f$ is called
\texttt{Walsh polynomial }if only  finitely many coefficients in its
Walsh representation are non-zero.

 The Walsh system is a dyadic analog of the trigonometric system, in particular,
\begin{equation}
\label{hat_jk}
	 \widehat{f_{j,n}}(\xi)= 2^{-j/2} {\rm w}_{n}(D^{-j} \xi)
	 \widehat{f}(D^{-j} \xi), \quad n\in \z_+, j\in\z.
\end{equation}

The function
$$
f^{[1]}(x):=\lim_{n\to \infty} \sum_{j=-n}^{n} 2^{j-1}(f(x)-f_{0, 2^{-j-1}}(x))
$$
is called \texttt{the Gibbs derivative} of $f$.
The inherited properties are
$$
 (f+g)^{[1]}=f^{[1]}+g^{[1]},\quad (cf)^{[1]}=c f^{[1]}, \quad
\widehat{f^{[1]}}\,=\,\lambda \widehat{f}
$$
\begin{equation}
{\rm w}^{[1]}_n\,=\, n {\rm w}_n,
\label{01}
\end{equation}

Unfortunately the Gibbs derivative does not inherit some natural
properties such as the chain rule and the rule $(f g)'\,=\,fg'+f'g.$
Moreover, Gibbs differentiability is not a local property of a function.
Higher order derivatives are defined recursively: $f^{[n+1]}=(f^{[n]})^{[1]}.$

\subsection{Haar system}
\label{Haar}

Lang~\cite{Lang} constructed
the Haar basis for $L_2(G)$
as one of the MRA-based wavelet systems on $G$.
Not only the basis itself (which also appeared
in the literature before~\cite{Lang}),
but some details of the construction will be
useful for us.

Set $\varphi=\mathbbm{1}_I$. This function has two wonderful
(from the point of view of wavelet theory) properties:
1) the functions $\varphi_{0k}$, $k\in\z_+$, form an orthonormal system;
2) $\varphi$ is a solution of the following refinement equation
\begin{equation}
\label{02}
2^{1/2} \varphi(x)=\varphi_{1,0}(x)+\varphi_{1,1}(x), \quad x\in G.
\end{equation}
It follows that $\varphi$ generates a multiresolution analysis  (MRA)
$\{V_j\}_{j\in\z}$, where
$$
V_j=\overline{\mbox{span}\{\varphi_{j,k}:\ k\in\z_+\}},\quad j\in\z.
$$
The union of all  spaces $V_j$ is dense in ${ L}_2(G)$, the intersection
of  all  spaces $V_j$ consists only of zero function, and
$V_j\subset V_{j+1}$ for all $j$.

The wavelet spaces $W_j$ are defined by
$$
 W_j= V_{j+1}\ominus V_j,\quad m\in \z,
$$
which yields the following orthogonal decomposition
\begin{equation}
\label{3}
L_2(G)={\bigoplus\limits_{j\in\z}W_j}.
\end{equation}
The corresponding wavelet function $\psi$ is given by
\begin{equation}
\label{03}
2^{1/2} \psi(x):=\varphi_{1,0}(x)-\varphi_{1,1}(x), \quad x\in G.
\end{equation}
The functions $\psi_{j,k}$, $j\in\z$, $k\in\z_+$,
form an orthonormal basis (Haar basis)
for $L_2(G)$, and any function  $f \in L_2(G)$ can be decomposed as
$$
f\ravno^{L_2}\sum_{j,k} a_{j,k}\psi_{j,k}
$$
(\texttt{Haar representation }in the sequel). A function $f$ is called
\texttt{Haar polynomial }if only finitely many coefficients in its
Haar representation are non-zero.

The functions  $\psi_{j,k}$, $j\in\z\setminus\z_+$, $k\in\z_+$, form an
orthonormal basis for $V_0$. Replacing these functions by the functions
$\varphi_{0,k}$, $k\in\z_+$, which also form an orthonormal basis for $V_0$,
we have another representation for $f \in L_2(G)$:
$$
f\ravno^{L_2}\sum_{k\in\mathbb{Z}_+} a_{k}\varphi_{0,k}+
\sum_{j,k\in\mathbb{Z}_+} a_{j,k}\psi_{j,k}
$$
(\texttt{quasi-Haar representation} in the sequel).

The MRA we described is an analog of the Haar MRA in the real setting,
but it is also an analog of the real Shannon MRA because
$\widehat\varphi=\varphi$. So, the space $V_j$ consists
of all functions $f\in L_2(G)$
whose Fourier-Walsh transform $\widehat f$ is supported in $I_{-j},$
the space $W_j$ consists of the functions whose Fourier-Walsh transform is supported in  $I_{-j-1}\setminus I_{-j},$ in particular,
\begin{equation}
\label{04}
\supp\widehat{\varphi_{j,k}}\subset I_{-j}, \quad k\in\z_+,
\end{equation}
\begin{equation}
\label{05}
\supp\widehat{\psi_{j,k}} \subset I_{-j-1}\setminus I_{-j}, \quad k\in\z_+.
\end{equation}
Moreover, taking into account~(\ref{hat_jk}), we have
\begin{equation}
\label{hat_Haar}
	\widehat{\varphi_{j,k}}(\xi)=2^{-j/2}{\rm w}_k (D^{-j}\xi)
\mathbbm{1}_{I_{-j}}(\xi)
	\quad
	\widehat{\psi_{j,k}}(\xi)=2^{-j/2}{\rm w}_k (D^{-j}\xi)
\mathbbm{1}_{I_{-j-1}\setminus I_{-j}}(\xi), \quad \xi\in G.
\end{equation}

\section{Distributions and functional classes on the Cantor group}
\label{classes}

We define  distributions on $G$, which can be considered as an
analog of tempered distributions in real analysis. Such a class
of distributions was introduced in two different ways
in the books~\cite{G} and~\cite{A}. We follow the latter one.

A function $\phi:\  G\to\cn$ is said to be locally constant
if for every $x\in G$ there exists a ball  containing $x$ such that $\phi$
is constant on this ball.
Any locally constant function is continuous on $G$.
A function $\phi$ is uniformly locally constant with rang $n$ if
it is constant on any ball of radius $2^{n}$.
We denote by $S$  the class of locally constant compactly supported  functions.
Evidently, any function in $S$ is uniformly locally constant.

\begin{prop}\cite[ \S 6.2]{GES} Let a function
  $\phi\in L_1(G)$ be  continuous.  Then $\phi$
  is uniformly locally constant with rang $n$
 if and only if  $\widehat\phi$ is supported in
 $I_n(0)$; and $\phi$
  is  supported  in $I_n(0)$ if and only if
  $\widehat\phi$  is uniformly locally constant with rang $n$.
\label{p1}
 \end{prop}
\begin{cor}
 A continuous function $\phi\in L_1(G)$  is in
$S$ if and only if  $\widehat{\phi}$ is in $S$.
\label{c1}
\end{cor}
The statement of Corollary follows immediately from Proposition~\ref{p1}.
Thus  $S$ is invariant with respect to the Walsh-Fourier transform,
and  it  is natural to say that $S$ is an analog of
the Schwartz class in the real analysis.

The convergence in $S$ is defined as follows. One says that a sequence
$\{\phi_k\}_{k=1}^\infty\subset S$ converges  to zero in $S$ if
\begin{enumerate}
\item
there exists a compact set $K\subset G$ such that $\supp \phi_k\subset K$
for all $k\in\n$;
\item
there exists $n\in\z$ such that any $\phi_k$ is constant on any ball of radius $2^n$;
\item
$\phi_k$ uniformly tends to zero on $K$ as $k\to\infty$.
\end{enumerate}
Evidently, the space $S$ is complete.

Let $S^\prime$ denote the dual space for $S$, i.e. $S^\prime$ consists of
continuous linear functionals $f:\ \phi\to \langle f, \phi\rangle$ on $S$.
The elements of  $S^\prime$ are called tempered distributions.
As usual, the convergence in $S^\prime$ is defined as the week convergence, i.e.
a sequence $\{f_k\}_{k=1}^\infty$ converges to zero in $S^\prime$, if
$$
\langle f_k, \phi\rangle\too\limits_{k\to\infty} 0\quad\forall\ \phi\in S.
$$
The completeness of $S^\prime$ can be checked in the standard way.

To give a representation for the tempered distributions
 we need the following axillary statement.

\begin{prop}
\label{charS}
If a function $\phi$ is in $S$, then its quasi-Haar representation is finite,
i.e. only a finitely many its coefficients
$\langle \phi, \psi_{jk}\rangle$ and $\langle \phi, \varphi_{0k}\rangle$
are non-zero. Conversely, any function $\phi$
with a finite quasi-Haar representation is in $S$.
\end{prop}
\textbf{Proof.}
If $\phi\in S$, then  $\phi\in L_2(G)$ and its Fourier-Walsh
transform $\widehat\phi$ is compactly supported due to Corollary~\ref{c1}.
It follows from the Plancherel theorem and~(\ref{05}) that
$$
\langle \phi, \psi_{jk}\rangle=\langle \widehat\phi, \widehat\psi_{jk}\rangle=0
$$
for all large enough $j$ and all $k\in\z_+$. On the other hand, for every $j$
there exist only finitely many $k\in\z_+$ such that
$\langle \phi, \psi_{jk}\rangle\ne0$ because  $\supp \phi$ is disjoint to
$\supp \psi_{jk}$ for large enough~$k$. Similarly, there exist
only finitely many $k\in\z_+$ such that $\langle \phi, \varphi_{0k}\rangle\ne0$.
Hence
$$
\phi=\sum_{k\in\mathbb{Z}_+} \langle \phi, \varphi_{0k}\rangle\varphi_{0k}+
\sum_{j,k\in\mathbb{Z}_+} \langle \phi, \psi_{jk}\rangle\psi_{jk},
$$
where both  sums are finite,
what proves the first statement. The second statement is trivial because
each of the functions $\varphi_{0k}, \psi_{jk}$ is in $S$ by definition.
$\Box$

The Fourier-Walsh transform $\widehat f$ of $f\in S'$ is defined by
$$
\langle \widehat f,\phi\rangle=\langle  f, \widehat\phi\rangle,\quad \phi\in S.
$$
The inversion formula~(\ref{12}) holds true  for any $f\in S'$ because the same is valid for any function in $S$.

If now $f$ is in $S^\prime$, $\phi\in S$, then, by Proposition~\ref{charS},
$$
\langle f, \phi\rangle=\sum_{k\in\mathbb{Z}_+}
\langle \varphi_{0k}, \phi\rangle\langle f, \varphi_{0k}\rangle+
\sum_{j,k\in\mathbb{Z}_+} \langle  \psi_{jk}, \phi\rangle\langle f, \psi_{jk}\rangle,
$$
where each of the sums is finite. Hence $f$ can be identified with the formal series
\begin{equation}
\sum_{k\in\mathbb{Z}_+} a_{k}\varphi_{0k}+
\sum_{j,k\in\mathbb{Z}_+} a_{j,k}\psi_{jk},\quad
a_{k}=\langle f, \varphi_{0k}\rangle,\quad a_{jk}=\langle f, \psi_{jk}\rangle,
\nonumber
\end{equation}
which is called     \texttt{quasi-Haar representation} of $f$.
 For convenience we will write
 $f=\sum_{k\in\mathbb{Z}_+} a_{k}\varphi_{0k}+
\sum_{j,k\in\mathbb{Z}_+} a_{j,k}\psi_{jk}$, and even
\begin{equation}
f(x)=\sum_{k\in\mathbb{Z}_+} a_{k}\varphi_{0,k}(x)+
\sum_{j,k\in\mathbb{Z}_+} a_{j,k}\psi_{j,k}(x),\quad  a_{k}, a_{j,k}\in\cn,\quad x\in G.
\label{1111}
\end{equation}

The space $S$ has "sufficiently many" functions to allow us to "tell regularly distributions apart". More precisely, the following analog of the du Bois-Reymond lemma  holds.
\begin{prop}
\label{p3}
 If a function $f$ is locally integrable on $G$ and its  quasi-Haar representation is zero, then $f=0$  almost everywhere on $G$.
\end{prop}
\textbf{Proof.}
Set $n\in\z_+$, $J_n=I_0(\lambda^{-1}(n))$, $g=f\varphi_{0n}$. Since
$\supp g\subset J_n$, $f=g$ on  $J_n$,
and each of the functions $\varphi_{0,k}$, $\psi_{j,k}$, $k,j\in\z_+$,  either
is supported or vanishes identically on $J_n$, the quasi-Haar representation of $g$ is also zero. On the other hand (see~\cite[Sec. 10.3.1]{GES}),
$$
g\ravno^{L(J_n)}\sum_{k\in\mathbb{Z}_+} \langle g, \varphi_{0k}\rangle\varphi_{0,k}+
\sum_{j,k\in\mathbb{Z}_+} \langle g, \psi_{jk}\rangle\psi_{j,k}.
$$
Hence $f=g=0$  almost everywhere on $J_n$, what was to be proved.
$\Box$

Denote by $\tilde{S}$ the set of Haar polynomials.
 Any Haar polynomial is in $L_2(G)$ and
 compactly supported. Hence its quasi-Haar representation
is also finite. It follows from Proposition~\ref{charS} that
$\tilde{S}\subset S$, i.e. $\tilde{S}$ is a subspace of $S$.
Denote by ${\tilde S}^\prime$ the space of continuous linear
functionals on $\tilde{S}$.

Since any $\phi\in\tilde{S}$ is a finite linear combination of the functions
$\psi_{jk}$,  $j\in \z$, $k\in\z_+$,
to define $f\in{\tilde{S}}^\prime$  on $\tilde{S}$ it suffices to define
$\langle f, \psi_{jk}\rangle$ for all $j\in\z$, $k\in\z_+$.
 Hence we can identify $f$ with  the formal series
\begin{equation}
\sum_{j,k\in\mathbb{Z}_+} a_{j,k}\psi_{j,k},\quad
 a_{j,k}=\langle f, \psi_{jk}\rangle,
\nonumber
\end{equation}
which is called \texttt{ Haar representation}  of  $f$.
Again  we can  write for convenience
\begin{equation}
\label{Hrep}
f(x)=\sum_{j,k\in\mathbb{Z}_+} a_{j,k}\psi_{j,k}(x),
\quad a_{j,k}\in\cn,\quad x\in G.
\end{equation}

Evidently, if $f\in S^\prime$, then $f\big|_{\tilde{S}}$ is in ${\tilde S}^\prime$.
But it follows from the next theorem that every element of ${\tilde S}^\prime$ is a
restriction of some $f\in S^\prime$ to ${\tilde{S}}$, i.e.
$$
{\tilde S}^\prime=\{f\big|_{\tilde{S}}\ , f\in S^\prime\}.
$$

\begin{teo}
\label{oneparam}
If $f\in \tilde{S}',$ then there exists a one-parametric family $\{f_c\}_{c\in\mathbb{C}}\subset S'$ such that $f_c\left|_{\tilde{S}}=f\right.$
for any $c\in \mathbb{C}$, and $g\left|_{\tilde{S}}\ne f\right.$ for any
$g\in S'\setminus \{f_c:\ c\in\mathbb{C}\}$.
Moreover, if (\ref{Hrep}) is the Haar representation of $f,$
then the quasi-Haar representation of $f_c$ is given by
$$
f_c(x)= \sum_{k=0}^{\infty} a_k(c)\varphi_{0,k}(x)+
\sum_{j,k \in \mathbb{Z}_+} a_{j,k} \psi_{j,k}(x), \ x\in G,
$$
where  $a_k=a_k(c)=a_k(0)+c,$ $k \in \z_+,$ $c\in \cn$,  is the general solution of system (\ref{syst}).\end{teo}

\textbf{Proof.}
First, we need to find a connection between the functions $\psi_{j,k}$, $j<0$, and $\varphi_{0,k}$.
It follows from (\ref{02}) and (\ref{03}) that
$$
2^{1/2}\psi_{j-1,k}=\varphi_{j,2k}-\varphi_{j,2k+1}, \quad
2^{1/2}\varphi_{j-1,k}=\varphi_{j,2k}+\varphi_{j,2k+1}.
$$
Iterating these recursion formulas,  we have
$$
2^{j/2}\psi_{-j,k}=\sum_{m=0}^{2^{j-1}-1}\varphi_{0,2^j k+m}-
\sum_{m=2^{j-1}}^{2^{j}-1}\varphi_{0,2^j k+m},
\quad j\in \mathbb{N},  \quad k\in \mathbb{Z}_+.
$$
Multiplying this by  a function $f_c\in S^\prime$, we obtain
$$
2^{j/2}\langle f_c,\psi_{-j,k}\rangle=
\sum_{m=0}^{2^{j-1}-1}\langle f_c, \varphi_{0,2^j k+m}\rangle-
\sum_{m=2^{j-1}}^{2^{j}-1}\langle f_c, \varphi_{0,2^j k+m}\rangle,
\quad j\in \mathbb{N},  \quad k\in \mathbb{Z}_+.
$$
The equality $f_c\left|_{\tilde{S}}=f\right.$ is equivalent to
$\langle f_c, \psi_{j,k}\rangle=\langle f, \psi_{j,k}\rangle =a_{jk}$,
 $j\in\z$, $k\in\z_+$.
 So, to find $f_c$ we have to solve
  the following  system with unknowns $a_{l},$ $l \in \mathbb{Z}_+$,
\begin{equation}
\label{syst}
	\sum_{m=0}^{2^{j-1}-1}a_{2^j k+m}-
	\sum_{m=2^{j-1}}^{2^{j}-1}a_{2^j k+m} = 2^{j/2} a_{-j,k}, \quad j\in\n, \ k\in\z_+,
\end{equation}
and set $\langle f_c, \varphi_{0,l}\rangle=a_l$, $l\in\z_+$,
$\langle f_c, \psi_{j,k}\rangle=a_{jk}$, $j,k\in\z_+$.

Choose an arbitrary parameter $c\in \cn$ and set $a_0=c$.
Then substituting $(j,\,k)=(1,\,0)$
in (\ref{syst}), we get
$
2^{1/2}a_{-1,0}=a_0-a_1.
$
So, $a_1 = a_0 - 2^{1/2} a_{-1,0}.$
On the second step we find the coefficients $a_2, a_3$ using (\ref{syst}) for
$(j,\,k)=(1,\,1)$ and $(j,\,k)=(2,\,0)$
$$
\left(
\begin{array}{cc}
1 & -1  \\
-1 & -1 \\
\end{array}
\right)
\left(
\begin{array}{c}
a_{2}  \\
a_{3} \\
\end{array}
\right)=
\left(
\begin{array}{c}
2^{1/2} a_{-1,1}  \\
2 a_{-2,0}-a_0-a_1 \\
\end{array}
\right)
$$
On the $J$-th step, $J\in\n$, the coefficients $a_k,$ $k=0,\dots,2^{J-1}-1$, are
known, and we find the coefficients $a_l,$ $l=2^{J-1},\dots,2^J-1$,
using (\ref{syst}) for
$(j,\,k)$ satisfying
\begin{equation}
\label{07}
2^{J-1}<2^j(k+1)\leq 2^J,\quad j\in\n, \ k\in\z_+.
\end{equation}
The number of equations, that is the number of all  solutions of~(\ref{07}), equals  the number of unknown coefficients that is $2^{J-1}.$

Consider the $J+1$-th step, $J\in\n.$
It follows from (\ref{syst}) that
$$
	\sum_{m=0}^{2^{j-1}-1}a_{2^j k+m+2^{J-1}}-
	\sum_{m=2^{j-1}}^{2^{j}-1}a_{2^j k+m+2^{J-1}} =2^{j/2}a_{-j,k+2^{J-1-j}},
$$
$$
	\sum_{m=0}^{2^{j-1}-1}a_{2^j k+m+2^{J}}-
	\sum_{m=2^{j-1}}^{2^{j}-1}a_{2^j k+m+2^{J}} =2^{j/2}a_{-j,k+2^{J-j}},
$$
where $j\in\n, \
	k\in\z_+,\ 2^{J-1} < 2^j(k+1)\leq 2^J, \ 1\leq j\leq J-1.$
Therefore, on $J+1$-th step, all but two equations for the unknowns
$a_{n_1},$ $n_1=2^{J},\dots,2^J+ 2^{J-1}-1$ and
$a_{n_2},$ $n_2=2^J+ 2^{J-1},\dots,2^{J+1}-1$
can be written using the equations of the previous $J$-th step written for the unknowns  $a_r,$ $r=2^{J-1},\dots,2^{J}-1,$ where $n_1=r+2^{J-1},$ and $n_2=r+2^{J}.$
To do so, it is sufficient to replace the right-hand side $2^{j/2}a_{-j,k}$ with $2^{j/2}a_{-j,k+2^{J-1-j}}$ and with
$2^{j/2}a_{-j,k+2^{J-j}}$ respectively.
The two remaining equations of the $J+1$-th step correspond to $(j,\,k)=(J,\,1)$ and $(j,\,k)=(J+1,\,0)$ and have the form
$$
	\sum_{m=0}^{2^{J-1}-1}a_{2^J +m}-
	\sum_{m=2^{J-1}}^{2^{J}-1}a_{2^J +m} =2^{J/2}a_{-J,1},
	$$
	$$
	-\sum_{m=2^{J}}^{2^{J+1}-1}a_{m}
	=2^{(J+1)/2}a_{-J-1,0}-\sum_{m=0}^{2^{J}-1}a_{m}.
$$

Denote by $M_J$  the matrix of the system  for  the $J$-th step.
By observation written in the previous paragraph, we see that the matrix $M_{J+1}$ is organized as follows.
We take the first  $2^{J-1}-1$ rows of the matrix $M_{J}$ and extend each row with $0$'s. The next $2^{J-1}-1$ rows are started with $0$'s and continued with the first $2^{J-1}-1$ rows of the matrix $M_J$. The first $2^{J-1}$ entries of the next  row equal $1$, the last $2^{J-1}$ entries of this row equal $-1$. Finally, the last row consists of $-1$'s. For example,
$$
M_3=\left(
\begin{array}{cccc}
 1 & -1 & 0 & 0 \\
 0 & 0 & 1 & -1 \\
1 & 1 & -1 & -1 \\
-1 & -1 & -1 & -1 \\
\end{array}
\right)
\quad
M_4=
\left(
\begin{array}{cccccccc}
1 & -1 & 0 & 0 & 0 & 0 & 0 & 0\\
 0 & 0 & 1 & -1 & 0 & 0 & 0 & 0  \\
1 & 1 & -1 & -1 & 0 & 0 & 0 & 0 \\
0 & 0 & 0 & 0 & 1 & -1 & 0 & 0 \\
0 & 0 & 0 & 0 & 0 & 0 & 1 & -1 \\
0 & 0 & 0 & 0 & 1 & 1 & -1 & -1 \\
1 & 1 & 1 & 1 & -1 & -1 & -1 & -1 \\
-1 & -1 & -1 & -1 & -1 & -1 & -1 & -1 \\
\end{array}
\right)
$$
To verify that the determinant of the matrix $M_J$ is not equal to zero we prove that $\det M_{J+1}=2(\det M_J)^2$ for $J>1.$
 Indeed,
let us consider the matrix $M_{J+1}$.
If we subtract the last row from the next to last one, divide the difference by $-2$, and insert the result
between the $2^J-1$ and the $2^J$ rows, we obtain a matrix that consists of four blocks. The blocks lying on the main diagonal are $M_J$ and one of the remaining blocks consists of $0$'s. Therefore, the determinant of this matrix is
$(\det M_J)^2.$  Since $\det M_1=-1,$ $\det M_2=-2,$ we obtain
$\det M_J=2^{2^{J-1}-1}$ for $J>2.$

Illustrate this  for the case  $J=2:$
$$
\det M_3 =
\left|
\begin{array}{cccc}
 1 & -1  & 0 & 0 \\
 0 & 0 & 1 & -1 \\
2 & 2  & 0 & 0 \\
-1 & -1 & -1 & -1 \\
\end{array}
\right|
=
-2 (-1)
\left|
\begin{array}{cccc}
 1 & -1 & 0 & 0 \\
 -1 & -1 & 0 & 0 \\
 0 & 0 & 1 & -1 \\
-1 & -1 & -1 & -1 \\
\end{array}
\right|
=
2
\left|
\begin{array}{cc}
 M_2 & 0 \\
 \dots & M_2 \\
\end{array}
\right|= 2 (\det M_2)^2.
$$

Thus, if we fix a parameter $c\in \cn$, then  system~(\ref{syst})
 has a unique solution $a_k$, $k\in \z_+$, depending on $c$. Denote by $a_k(0)$
the solution corresponding to $c=0$ and prove
 that the general solution is given by
$a_0(c)=c,$ $a_k(c)=a_k^0+c$, $k\in\n,$ $c\in \cn.$
Set  $b_k:=a_k(c)-a_k(0)$. Due to~(\ref{syst}),
we have
\begin{equation}
\label{08}
B:=	\sum_{m=0}^{2^{j-1}-1}b_{2^jk+m}-
	\sum_{m=2^{j-1}}^{2^{j}-1}b_{2^jk+m} =0, \quad j\in\n, k\in\z_+.
\end{equation}
Let us prove  by induction on $j$  that
\begin{equation}
\label{09}
b_{2^jk}=b_{2^jk+1}=\dots=b_{2^jk+2^j-1}
 \end{equation}
for every $k\in\z_+$. The induction base for $j=1$ follows from~(\ref{08}) immediately.
Let us check the induction step from $j-1$ to $j$. Due to the induction hypothesis,
$$
B=\sum_{m=0}^{2^{j-1}-1}b_{2^{j-1}(2k)+m}-
	\sum_{m=2^{j-1}}^{2^{j}-1}b_{2^{j-1}(2k+1)+m}=
2^{j-1}(b_{2^{j-1}(2k)}-b_{2^{j-1}(2k+1)}).
$$
This and~(\ref{08}) yield  $b_{2^{j-1}(2k)}=b_{2^{j-1}(2k+1)}$,
which together with the induction hypothesis proves~(\ref{09}).
It follows that $b_l=b_0=c$ for every $l\in\n$. Hence
the quasi-Haar representation of $f_c-f_{0}$ is $\sum_{k=0}^{\infty} c\varphi_{0,k}$, what was to be proved. $\Box$

We next introduce a class of periodic distributions.
Denote by $P$  the class of $1$-periodic locally constant
functions. It is clear that any $\phi\in P$ is uniformly locally constant.
The convergence in $P$ is defined as follows. One says that a sequence
$\{\phi_k\}_{k=1}^\infty\subset P$ converges  to zero in $P$ if
\begin{enumerate}
\item
there exists $n\in\z$ such that any $\phi_k$ is constant on any ball of radius $2^n$;
\item
$\phi_k$ uniformly tends to zero on $I$ as $k\to\infty$.
\end{enumerate}
Evidently, the space $P$ is complete.

Let $P^\prime$ denote the dual space for $P$, i.e. $P^\prime$ consists of
continuous linear functionals $f:\ \phi\to \langle f, \phi\rangle$ on $P$.
The convergence in $P^\prime$ is defined as the week convergence, i.e.
a sequence $\{f_k\}_{k=1}^\infty$ converges to zero in $P^\prime$, if
$$
\langle f_k, \phi\rangle\too\limits_{k\to\infty} 0\quad\forall\ \phi\in P.
$$
The completeness of $P^\prime$ can be checked in the standard way.

To give a representation for these distributions
 we need the following simple statement.

\begin{prop}\cite[ \S 1.4, \S 2.7]{GES}
\label{p2}
A $1$-periodic function $\phi$ is in $ P$ if and only if it is a Walsh polynomial.
\end{prop}

Since, by Proposition~\ref{p2}, any $\phi\in P$ is a finite linear combination
of the functions ${\rm w}_k$,  $k\in\z_+$, to define $P^\prime$  on $P$
it suffices to define  $\langle f, {\rm w}_{k}\rangle$ for all  $k\in\z_+$.
 Hence we can identify $f$ with  the formal series
\begin{equation}
\sum_{k=0}^{\infty} a_{k}{\rm w}_{k},\quad
 a_{k}=\langle f, {\rm w}_{k}\rangle,
\nonumber
\end{equation}
which is called its\texttt{ Walsh representation}.
For convenience we  write
\begin{equation}
\label{fWs}
f(x)=\sum_{k=0}^{\infty} c_{k}{\rm w}_{k}(x),\quad c_k\in\cn, \quad x\in G.
\end{equation}

Finally, we introduce an operation of differentiation for periodic distributions.
If $\phi$ is a Walsh polynomial, then its Gibbs derivative $\phi^{[1]}$ is also
a Walsh polynomial, and we can well-define the Gibbs derivatives of
$f\in P^\prime$ by
$$
\langle f^{[k]}, \phi\rangle=\langle f, \phi^{[k]}\rangle,\quad \phi\in P.
$$

 \section{Gibbs differential equations}
\label{Gde}
 Consider a linear  differential equation of order $n$
 \begin{equation}
\label{ode_n}
	\sum_{k=0}^{n}\alpha_k f^{[k]}\,=\,g, \quad n \in \mathbb{N},
\quad \alpha_{k} \in \mathbb{C}, k=0,\dots, n.
\end{equation}
As usual, the characteristic polynomial associated with this equation is
$Q_n(y)=\sum_{k=0}^{n}\alpha_k y^k$.
 \begin{teo}
\label{odeN}
Let  $g$ be a distribution in $P'$,
$\sum_{m=0}^{\infty} d_m  {\rm w}_m$, be
the Walsh representation of $g$, and let $Q_n$ be the characteristic polynomial
associated with~(\ref{ode_n}). Then
\begin{enumerate}
	\item
	if $Q_n$ has no integer non-negative roots, then (\ref{ode_n})
	has a unique solution $f=\sum_{m=0}^{\infty} c_m {\rm w}_m$ in~$P',$
	where $c_m=d_m/Q_n(m)$;
	\item
	if $\{m_1,\dots, m_{r}\}$  is the set of integer non-negative roots of $Q_n$
	  and $\sum_{l=1}^r |d_{m_l}| = 0,$
	 then (\ref{ode_n})
	has an $r$-parametric family (enumerated by $c_{m_1},\,\dots,\,c_{m_r}$) of solutions  in~$P'$
	$$
     f(\cdot,\, c_{m_1},\,\dots,\,c_{m_r})=\sum_{m\in\z_+} c_m {\rm w}_m,
		$$
	where $c_m=d_m/Q_n(m)$ for $m\neq m_l,$ $l=1,\dots,r$, and $c_{m_1},\,\dots,\,c_{m_r}$ are arbitrary complex numbers;
	\item
if $\{m_1,\dots, m_{r}\}$  is the set of integer non-negative roots of $Q_n$
	  and $\sum_{l=1}^r |d_{m_l}|\neq 0,$
	 then (\ref{ode_n})
	has no solutions  in $P'$.	
\end{enumerate}
Moreover, all solutions from items $1$ and $2$ are in $L_2^P$ whenever
\begin{equation}
	\label{odeL2}
	d_m = O(m^{n-1/2-\varepsilon}), \quad \varepsilon>0,\quad  m\to \infty;
\end{equation}
	solutions all from items $1$ and $2$ are in $C^P$ whenever
	\begin{equation}
	\label{odeC}
		d_m = O(m^{n-1-\varepsilon}), \quad \varepsilon>0,\quad   m\to \infty.
\end{equation}

\end{teo}

\textbf{Proof.}
If $\sum_{m=0}^{\infty}c_m {\rm w}_m$ is the Walsh representation
of $f \in P',$ then, taking into account that,
by (\ref{01}),
${\rm w}_m^{[k]}=m^k {\rm w}_m$, we can
 rewrite~(\ref{ode_n}) in the form
$$
\sum_{k=0}^{n}\alpha_k \sum_{m=0}^{\infty}c_m m^k {\rm w}_m=
\sum_{m=0}^{\infty}d_m  {\rm w}_m,
$$
or equivalently
$$
\sum_{m=0}^{\infty}c_m Q_n(m) {\rm w}_m=
\sum_{m=0}^{\infty}d_m  {\rm w}_m.
$$
Therefore, we get the equations in $c_m$
$$
c_m Q_n(m) =d_m \quad m\in\z_+.
$$
Proof of items $1$-$3$ follows immediately. Also it follows that
$$
c_m=d_m/Q_n(m)=O( d_m m^{-n}),\quad m\to\infty.
$$
So, if~(\ref{odeL2}) is fulfilled, then
$\sum_{m=0}^{\infty}|c_m|^2<\infty$, which yields $f\in L_2^P$.
If~(\ref{odeC}) is fulfilled, then
$\sum_{m=0}^{\infty}|c_m|<\infty$, which yields the uniform convergence
of the series $\sum_{m=0}^{\infty} c_m {\rm w}_m$, and its sum is continuous.
\hfill$\Box$

  \begin{rem}
   If in~(\ref{ode_n}) the function $g$ is in $L_1^P$, then $d_m=o(1)$, as $m\to\infty$. It follows that~(\ref{odeL2}) is fulfilled for any $n\in \n$, and~(\ref{odeC}) is fulfilled whenever $n>1$. Hence,  in this case all periodic solutions of~(\ref{ode_n}) are   in $L_2^P$, and all periodic solutions are continuous functions whenever $n>1$.
  \end{rem}

We now consider the Cauchy problem  for the one-dimensional homogeneous wave equation in  variables $(x,t)$,  where $x\in G$ and  $t$ (time) is real. Let  $U=[0,\,+\infty)$  or  $U=[0,\,T].$
The problem is of the form
\begin{equation}
	\label{pdeGibbswave}
	\left\{
	\begin{array}{l}
	\frac{\partial ^2 f(x,\,t)}{\partial t^2} = f^{[2]}_{x^{2}} f(x,\,t), 	 \\
	f(x,\,0) = f^0(x),\ \
	f'_t(x,\,0) = f^1(x),
	\end{array}
	\right.\quad x\in G,\quad t\in U.
\end{equation}

\begin{teo}
\label{teopdeGibbswave}
Let $f^0,$ $f^1 \in P',$
$\sum_{n=0}^{\infty}p_n {\rm w}_n$ and $\sum_{n=0}^{\infty} q_n {\rm w}_n$
be the Walsh representation of $f^0$ and $f^1$ respectively.  Then
\begin{enumerate}
\item
the Cauchy problem (\ref{pdeGibbswave}) has a unique solution $f(x,\,t)$ which is in $P'$
for every $t\in U;$
 \item
 this solution is  in $C^P$ whenever
\begin{equation}
	\label{pdeL2C}
	p_n=O(e^{-n \theta(n)}),\ \ q_n=O(e^{-n \theta(n)}), \quad n\to \infty,
\end{equation}
where $\theta(n)\to \infty$ as $n \to \infty;$
 \item this solution is in $L_2^P$ whenever (\ref{pdeL2C}) holds true.
\end{enumerate}
\end{teo}

\textbf{Proof.}
1.
Fix $t\in U.$
Suppose
$f(\cdot,\,t) = \sum_{n=0}^{\infty}c_{n}(t){\rm w}_n$
is the Walsh representation of $f(\cdot,t)$.
 Then, taking into account that,
by (\ref{01}),
${\rm w}_n^{[k]}=n^k {\rm w}_n$,  we can rewrite~(\ref{pdeGibbswave})
in the form
$$
\left\{
	\begin{array}{l}
	 \sum\limits_{n=0}^{\infty} \ddot{c}_n(t){\rm w}_n = \sum\limits_{n=0}^{\infty} n^2 c_n(t){\rm w}_n, \\
	 \sum\limits_{n=0}^{\infty}c_{n}(0){\rm w}_n = \sum\limits_{n=0}^{\infty}p_{n}{\rm w}_n \ \
	 \sum\limits_{n=0}^{\infty}\dot{c}_{n}(0){\rm w}_n = \sum\limits_{n=0}^{\infty} q_{n}{\rm w}_n \\
	\end{array}
	\right.
$$
where  $\dot{c}_{n}$ and $\ddot{c}_{n}$ denotes the ordinary first and second derivatives of  $c_{n}$ with respect to $t$.
Therefore, for every $n\in \n$ we obtain the  Cauchy problem for a linear ordinary differential equation of the second order
$$
\left\{
	\begin{array}{l}
\ddot{c}_n(t) = n^2 c_n(t),
	\\
	c_n(0)=p_n,
	\ \
	\dot{c}_n(0)=q_n.
	\end{array}
	\right.
$$
 The solution to the problem is
 $$
 c_n(t)= p_n \cosh (nt) +\frac{q_n}{n} \sinh (nt).
 $$
 Thus the coefficients in the Walsh representation of $f$ are found.

2. If~(\ref{pdeL2C}) is fulfilled, then
$\sum_{n=0}^{\infty}|c_n|<\infty$, which yields $f\in C^P$.

3. If~(\ref{pdeL2C}) is fulfilled, then
$\sum_{n=0}^{\infty}|c_n|<\infty$, so,
$\sum_{n=0}^{\infty}|c_n|^2<\infty$,
 which yields $f\in L_2^P$.
\hfill $\Box$

 \section{Modified Gibbs derivatives}
\label{MGd}
 In the previous section we described a method for finding $1$-periodic
 solutions to linear differential equations. Let us analyze this method for
 a simple equation $f^{[1]}=g$, where $g$ is $1$-periodic, say $g={\rm w}_k$.
 We find the Fourier-Walsh coefficients of $f^{[1]}$, that are
 $\widehat{f^{[1]}}(n)=\delta_{n,k}$,
 then, using $\widehat{f^{[1]}}(n)=n \widehat{f}(n)$, find the Fourier-Walsh
 coefficients of $f$, which allows to restore~$f$, that is $f={\rm w}_k/k$.

 Let us try to  repeat our trick for a non-periodic function $g,$ using the
 Fourier-Walsh transform instead of the Fourier-Walsh
 coefficients, i.e. find $\widehat f$ and then restore $f$.
 Suppose $g={\rm w}_k \varphi$, where, as above, $\varphi=\mathbbm{1}_I,$
 and try to find a  solution $f$ to the equation
 $f^{[1]}=g$. It would be desirable to obtain a solution supported on $I$.
 Since $\widehat{f^{[1]}}=\lambda \widehat{f},$ and
$\widehat{{\rm w}_k \varphi}=\varphi_{0,k},$ we have
$\widehat{f}=(1/\lambda)\varphi_{0,k}.$
Therefore,
$$
f(x)=\int\limits_{G}(\lambda(\xi))^{-1}\varphi_{0,k}(\xi){\rm w}(x,\,\xi)\,d\xi.
$$
The solution is not supported on $I$ because
 $$
 f(x)={\rm w}_k(x)\int_{I} (\lambda(\xi \oplus k))^{-1}{\rm w}_{\lfloor\lambda(x)\rfloor}(\xi)\,d\xi,
 $$
and $f(x)={\rm w}_k(x) \log ((k+1/2)^2 k^{-1}(k+1)^{-1})$ for
$x \in I_{-1}\setminus I$.  It is not difficult to check that
$f$ is even not a compactly supported functions. Moreover,
 $f(x)={\rm w}_k(x) \log ((k+1)/k)$ for  $x \in I,$
which yields that $f$  does not coincide on $I$
 with the periodic solution of the equation $f^{[1]}={\rm w}_k$
 However, if we  modify a little bit the definition of the derivative, then  the situation  cardinally changes. Let us use an operator ${\cal D}$  defined (on an appropriate class of functions) by
$$
{\cal D}f(x)=\int\limits_{G} \|\xi\| \widehat{f}(\xi) {\rm w}(\xi,\,x)\,d\xi
$$
instead of the Gibbs derivative.
Then the solution of the equation $Df={\rm w}_k \varphi$ is compactly supported. It takes the form
$$
f(x)=\int\limits_{G}\|\xi\|^{-1}\varphi_{0,k}(\xi){\rm w}(x,\,\xi)\,d\xi=\|\lambda^{-1}(k)\|^{-1}
\int\limits_{G}\varphi_{0,k}(\xi){\rm w}(x,\,\xi)\,d\xi=
\|\lambda^{-1}(k)\|^{-1}{\rm w}_k(x) \varphi(x).
$$

In a similar way one can introduce the modified derivative ${\cal D}$ for the $1$-periodic function.
 Indeed, suppose $f$ is a $1$-periodic function defined on an appropriate class,
 $\sum_{k \in \z_+} \widehat{f}(k) {\rm w}_k$ is its Walsh series; then
 $$
 {\cal D}f = \sum\limits_{k \in \z_+}\|\lambda^{-1}(k)\| \widehat{f}(k) {\rm w}_k
 $$
 It is clear that  ${\cal D} {\rm w}_k = \|\lambda^{-1}(k)\| {\rm w}_k,$ i.e.
 the Walsh functions are eigenfunctions of  ${\cal D}$ too.
  In contrast to the  Gibbs derivative,  if $g$ is a $1$-periodic function,  then the periodic solution $f$ of the equation ${\cal D}f = g$  coincides on $I$ with the solution of the equation ${\cal D}f = g \mathbbm{1}_I.$
 Indeed, if $\sum_{k \in \z_+} \widehat{g}(k) {\rm w}_k$ is the  Walsh series of $g$, then the equation
 ${\cal D}f = g$  can be rewritten as
 $$
 \sum\limits_{k \in \z_+} \|\lambda^{-1}(k)\| \widehat{f}(k) {\rm w}_k =
 \sum\limits_{k \in \z_+} \widehat{g}(k) {\rm w}_k,
 $$
 that is
 $
 \widehat{f}(k) = \|\lambda^{-1}(k)\|^{-1} \widehat{g}(k),
 $
 which yields
 $
 f = \sum_{k \in \z_+} \|\lambda^{-1}(k)\|^{-1} \widehat{g}(k) {\rm w}_k.
 $
 On the other hand,
 the equation
 ${\cal D}f = g \mathbbm{1}_I$
 can be rewritten as
 $
 \|\xi\| \widehat{f} = \widehat{g \mathbbm{1}_I},
 $
 that is
 $
 \|\xi\| \widehat{f} = \sum_{k\in \z_+}\widehat{g}(k) \varphi_{0,k}.
 $
 Hence
 $$
 f(x) = \int\limits_{G} \sum\limits_{k\in\z_+} \|\xi\|^{-1} \widehat{g}(k) \varphi_{0,k}(\xi) {\rm w}(\xi,\,x)\,d\xi=
 \sum\limits_{k\in\z_+} \widehat{g}(k) \int\limits_{G}  \|\xi\|^{-1}  \varphi_{0,k}(\xi) {\rm w}(\xi,\,x)\,d\xi
 $$
 $$
 =
 \sum\limits_{k\in\z_+} \widehat{g}(k) \int\limits_{I}  \|\xi \oplus \lambda^{-1}(k)\|^{-1}
 {\rm w}(\xi \oplus \lambda^{-1}(k),\,x)\,d\xi=
\sum\limits_{k\in\z_+} \widehat{g}(k) {\rm w}_k(x) \int\limits_{I}  \|\xi \oplus \lambda^{-1}(k)\|^{-1}
 {\rm w}(\xi,\,x)\,d\xi
 $$
  Taking into account that $\|\xi \oplus \lambda^{-1}(k)\| = \|\lambda^{-1}(k)\|$ for $\xi \in I$, we finally get
 $$
 f(x) = \sum\limits_{k\in\z_+} \widehat{g}(k) {\rm w}_k(x) \|\lambda^{-1}(k)\|^{-1} \int\limits_{I} {\rm w}(\xi,\,x)\,d\xi=
 \left(\sum\limits_{k\in\z_+}
  \|\lambda^{-1}(k)\|^{-1} \widehat{g}(k) {\rm w}_k(x)\right)
 \mathbbm{1}_I.
 $$
 Since this is just an illustration of the general idea we suppose that all operations with sums, integrals and so on are justified.

We now introduce fractional
modified Gibbs derivatives ${\cal D}^\alpha$ defined on ${\tilde S}^\prime$.
 This operator was introduced on $L_1(G)$  in~\cite{G}.
Such kind of operators are often called pseudo-differential.

First we need the following axillary statement.

\begin{lem}
\label{l1}
 Let $\phi\in S$. For $\phi$ to be in $\tilde S$ it is necessary and sufficient
 that
 \begin{equation}
  \supp\widehat \phi\cap I_n =\emptyset \quad\mbox{for some}\  n\in\z.
 \label{06}
 \end{equation}
\end{lem}
\textbf{Proof.} The necessity follows from~(\ref{05}).
Let~(\ref{06}) be satisfied for $\phi\in S$. It follows  from the Plancherel equality and ~(\ref{05})  that
$$
\langle \phi, \psi_{jk}\rangle=\langle \widehat\phi, \widehat\psi_{jk}\rangle=0,
\quad k\in\z_+,
$$
whenever  $-j$ is large enough. On the other hand, it was
shown in the proof of Proposition~\ref{charS} that $\langle \phi, \psi_{jk}\rangle=0$  for all $k\in\z_+$, whenever $j$ is large enough,
and for every $j$ there exist
only finitely many $k\in\z_+$ such that $\langle \phi, \psi_{jk}\rangle\ne0$.
Therefore,  $\phi$ is a Haar polynomial, which proves the sufficiency. $\Box$

Let  $\alpha\in\r$, ${\mathbb D}^\alpha(x):=\|x\|^\alpha$ for $x\in G$, $x\ne{\bf0}$,
and ${\mathbb D}^\alpha({\bf0}):=1$.
Define the \texttt{ fractional modified Gibbs derivative }${\cal D}^{\alpha}$
on $\tilde S$ by

$$
\widehat{{\cal D}^{\alpha}\phi}={\mathbb D}^\alpha\widehat{\phi}, \quad\phi\in\tilde S.
$$
Due to Lemma~\ref{l1}, ${\cal D}^\alpha$ is well defined and
${\cal D}^\alpha: \tilde S\to\tilde S$. Moreover, ${\cal D}^{-\alpha}$ is the inverse
operator to ${\cal D}^{\alpha}$, which yields that ${\cal D}^{\alpha}$ is a one-to-one
map taking $\tilde S$   onto $\tilde S$. This  allows to extend
 fractional modified Gibbs derivatives to $\tilde S^\prime$.
 For $f\in {\tilde S}^\prime$, we can well-define
 ${\cal D}^{\alpha}f\in {\tilde S}^\prime$ by
$$
\langle {\cal D}^{\alpha}f, \phi\rangle=\langle f, {\cal D}^{\alpha}\phi\rangle,\quad \phi\in \tilde S.
$$

  \begin{rem}
  Consider the  equation
\begin{equation}
	\label{A0}
	{\cal D}^0 f = g, \quad g \in \tilde{S}'.
\end{equation}
  Since ${\cal D}^0$ is the identical operator on $\tilde{S}',$ evidently, there is a unique solution $f=g$ in $\tilde{S}'.$ Any locally integrable function, in particular, a continuous one, is in $\tilde{S}'.$ Assume that $g$ is a continuous function. In this case it is natural to say that the solution $f$ is also a continuous function. However, by Theorem \ref{oneparam}, there exist infinitely many different continuous functions $f$ satisfying (\ref{A0}).
  All these functions have the same Haar representation, but they have different quasi-Haar representations. On the other hand,   if $g\in L_2(G)$, then only one  function $f \in L_2(G)$   satisfies (\ref{A0}).
  \end{rem}


\begin{prop}
\label{eigen}
Suppose $g,$ $\widehat{g},$ $\mathbb{D}^{\alpha} \widehat{g}$ are locally integrable on $G$, $j\in\z$. Then the assertion
${\rm supp}\, \widehat{g} \subset I_{-j-1}\setminus I_{-j}$
is necessary and
sufficient for
$g$ to be an eigenfunction of ${\cal D}^{\alpha}$
corresponding to the eigenvalue~$2^{j \alpha}$.
\end{prop}
\textbf{Proof.} Since the function $g$ is locally integrable, it is also in $\tilde{S}'$.
Assume that ${\rm supp}\, \widehat{g} \subset I_{-j-1}\setminus I_{-j}$.
Using the Plancherel  equality, for  any $\phi\in\tilde{S}$ we obtain
$$
\langle {\cal D}^{\alpha}g, \phi\rangle=\langle g, {\cal D}^{\alpha}\phi\rangle=
\langle \widehat g, \widehat{{\cal D}^{\alpha}\phi}\rangle=
\int\limits_{I_{-j-1}\setminus I_{-j}} {\widehat{g}(\xi)}{\|\xi\|^{\alpha}} \widehat{\phi}(\xi)\,d\xi=
2^{j\alpha}\int\limits_G {\widehat{g}(\xi)} \widehat{\phi}(\xi)\,d\xi=
2^{j\alpha}\langle g, \phi\rangle,
$$
which proves the necessity.

Let now $g$  be an eigenfunction of ${\cal D}^{\alpha}$
corresponding to the eigenvalue $2^{j_0 \alpha}$. Again by the Plancherel  equality, we have
$$
0=\langle {\cal D}^{\alpha}g, \phi\rangle-
2^{j\alpha}\langle g, \phi\rangle=
\int\limits_{G} {\widehat{g}(\xi)}({\|\xi\|^{\alpha}}-2^{j_0\alpha}) \widehat{\phi}(\xi)\,d\xi
$$
for  every $\phi\in\tilde{S}$. 
Therefore, for $\phi=\psi_{j,k}$ $j\in\z,$ $k\in\z_+$, we have
$$
0=\int\limits_{G} {\widehat{g}(\xi)}({\|\xi\|^{\alpha}}-2^{j_0\alpha}) \widehat{\psi}_{j,k}(\xi)\,d\xi= 2^{-j/2}
\int\limits_{G} {\widehat{g}(\xi)}({\|\xi\|^{\alpha}}-2^{j_0\alpha}) {\rm w}_k (D^{-j}\xi)
\mathbbm{1}_{I_{-j-1}\setminus I_{-j}}(\xi)\,d\xi
$$
$$
= 2^{-j/2}(2^{j\alpha}-2^{j_0\alpha})\int\limits_{I_{-j-1}\setminus I_{-j}}\widehat{g}(\xi)
{\rm w}_k (D^{-j}\xi)\,d\xi=
2^{j/2}(2^{j\alpha}-2^{j_0\alpha}) \int\limits_{I}\widehat{g}(2^j(\eta\oplus 1)) {\rm w}_k (\eta)\,d\eta.
$$
Since $\widehat{g}$ is locally integrable on $G$, 
the function $h(\eta)=\widehat{g}(2^j(\eta\oplus 1))$ is integrable on  $I$ 
for every fixed~$j$, and if $j\ne j_0$, then all Fourier-Walsh coefficients of $h$ 
equal zero. It follows that  $h=0$ almost everywhere on $I$ 
for $j\neq j_0$ (see~\cite[Sec. 10.2.1]{GES}).  Hence $\widehat{g}=0$
almost everywhere on $I_{-j-1}\setminus I_{-j}$ whenever 
$j\ne j_0$, which proves the sufficiency. 
\hfill $\Box$

\begin{rem}
If $g \in L_2(G)$, then the necessity of the conclusion of Proposition \ref{eigen} holds for any $\alpha \in \r$, and the sufficiency holds for  $\alpha>0$.
\end{rem}

\begin{cor}
\label{eigen_Haar}
Any Haar function $\psi_{j,k}$ is an eigenfunction of ${\cal D}^{\alpha}$ corresponding to the eigenvalue $2^{j \alpha}.$
\end{cor}

\textbf{Proof} The statement follows from Proposition \ref{eigen} and (\ref{hat_Haar}).\hfill $\Box$

\subsection{Modified Gibbs differential equations}

Consider the equation
\begin{equation}
	\label{odePDO}
	{\cal D}^{\alpha} f +\beta f = g, \quad \alpha, \beta\in\r
\end{equation}
with respect to the unknown function $f$.

\begin{teo}
\label{teo_odePDO}
Let $g$ be a distribution in $\tilde{S}'$,
 $\beta \neq -2^{j \alpha}$ for all $j \in \mathbb{Z}.$
Then
\begin{enumerate}
	\item equation (\ref{odePDO}) has a unique solution in $\tilde{S}'$;
    \item the solution is in $L_2(G)$ whenever $g\in L_2(G)$ and $\beta\ne0$;
	\item if $g$ is continuous on $G,$ and either $\beta\ne0$ or $\beta=0$,
     $\alpha>-1/2,$ 	then all continuous solutions of (\ref{odePDO}) are given by
\begin{equation}
\label{10}
f_c=\sum_{k=0}^{\infty}a_k\varphi_{0,k}+
\sum_{k,j\in \mathbb{Z}_+}
\frac{\left\langle g,\,\psi_{j,k}
\right\rangle}{2^{j \alpha}+\beta} \psi_{j,k}+c,
\end{equation}
where $c$ is an arbitrary complex number, and
$a_k$, $k\in\z_+$, is a solution of system (\ref{syst}).
All the functions $f_c$ have the same Haar representation.
\end{enumerate}
\end{teo}

\textbf{Proof.}
1.
Let $f=\sum_{j,k} a_{j,k}\psi_{j,k}$,
be the Haar representations of $f\in \tilde{S}'$.
It follows from Corollary~\ref{eigen_Haar} that equation (\ref{odePDO}) can be rewritten as
$$
\sum_{j,k}(2^{j \alpha} + \beta) a_{j,k} \psi_{j,k}=
\sum_{j,k}\left\langle g,\,\psi_{j,k}
\right\rangle \psi_{j,k},
$$
 or equivalently
\begin{equation}
\label{11}
(2^{j \alpha} + \beta) a_{j,k}=\left\langle g,\,\psi_{j,k}
\right\rangle, \quad j\in\z, \ \ k\in\z_+.
\end{equation}
It follows that
 $$
 f=\sum_{j\in\mathbb{Z},k\in\mathbb{Z}_+}
 \frac{\left\langle g,\,\psi_{j,k}
\right\rangle}{2^{j \alpha} + \beta}
 \psi_{j,k}.
  $$
  This proves item 1. If $\beta\ne0$, then $|{2^{j \alpha} + \beta}|\ge\delta$ for some $\delta>0$ and all $j\in\z$, which yields item~2.

  Let now the assumptions of item~3 be fulfilled. According to Theorem \ref{oneparam}, the functions $f_c$ given by~(\ref{10}),
and only these functions  in~$S'$, have the same Haar representation as $f$.
Let us prove that $f_0$ is continuous.

Fix $x_0\in G$ and the compact set $I_0(x_0).$
By the definition of the Haar function,
for every $x\in I_0(x_0)$
there is a unique $k=k(x)=k(x_0)\in \mathbb{Z}_+$ such that
$x\in {\rm supp}\,\varphi_{0,k(x_0)}$, and
for every $j \in \mathbb{Z}_+$  there is a unique $k=k(j,x)\in \mathbb{Z}_+$ such that $x\in {\rm supp}\,\psi_{j,k(j,x)}.$ So, using the quasi-Haar representation for $f_0$, we have
$$
f_0(x)=a_{k(x)} \varphi_{0,k(x)}(x) + \sum_{j\in \mathbb{Z}_+}
\frac{\left\langle g,\,\psi_{j,k(j,x)}
\right\rangle}{2^{j \alpha} + \beta}\psi_{j,k(j,x)}(x).
$$
Since $g$ is continuous,
$g$ is bounded on the compact set $I_0(x_0),$
that is $|g(x)|\leq M.$
Therefore,
$$
\left|\left\langle g,\,\psi_{j,k(j,x)}\right\rangle\right|=
\left|\int\limits_{G} g \psi_{j,k(j,x)} \right|\leq
M \int\limits_{G}  |\psi_{j,k(j,x)}| = M 2^{-j/2}.
$$
Thus
$$
\frac{|\left\langle g,\,\psi_{j,k(j,x)}\right\rangle|}{2^{j\alpha}+\beta}\leq M\frac{2^{-j/2}}{2^{j\alpha}+\beta}=
\left\{
	\begin{array}{l}
	O( 2^{j(-1/2-\alpha)}), 	\ \ \mbox{if}\ \ \ \alpha>0 \ or \ \  \beta=0,
\\
	O( 2^{-j/2}), \ \ \ \ \ \ 	\ \ \mbox{if}\ \ \ \alpha\le0, \  \beta\ne0,
	\end{array}
	\right.
 \mbox{ as } j \to \infty.
$$
Since the latter estimate is uniform on $I_0(x_0)$,  the series
$$
\sum_{j\in \mathbb{Z}_+}\frac{\left\langle g,\,\psi_{j,k(j,x)}\right\rangle}{2^{j \alpha}+\beta}
\psi_{j,k(j,x)}(x)
$$
uniformly converges on $I_0(x_0),$ which yields that $f_0$  is continuous on $I_0(x_0).$ It remains to note that $x_0$ is an arbitrary element of $G,$
thus $f_0$ is continuous on $G.$
\hfill  $\Box$

\begin{rem}
It is clear from the proof of Theorem~\ref{teo_odePDO}, that the assumption
of continuity of $g$ in item~3 can be replaced by the boundedness of $g$ on any compact set.
\end{rem}

\begin{teo}
Let $g$ be a distribution in $\tilde{S}'$,
  $\beta = -2^{j_0 \alpha},$ $j_0 \in \mathbb{Z}_+.$ If
 $\langle g,\,\psi_{j_0,k} \rangle \neq 0$ for at least one $k \in \mathbb{Z}_+,$ then (\ref{odePDO}) has no solutions in $\tilde{S}'.$
  If $\langle g,\,\psi_{j_0,k} \rangle = 0$ for all $k \in \mathbb{Z}_+,$ then
\begin{enumerate}
	\item equation (\ref{odePDO}) has a family of solutions in $\tilde{S}'$
$$
f_{\{c_k\}}=
\sum_{k\in\z_+, j\in \mathbb{Z}, j\neq j_0}
\frac{\left\langle g,\,\psi_{j,k}\right\rangle}{2^{j \alpha}+\beta} \psi_{j,k}
+\sum_{k=0}^{\infty}c_k\psi_{j_0,k},
$$
where $c_k$, $k\in\z_+$, are arbitrary complex numbers.
	\item if $g$ is continuous on $G,$ 	then all continuous solutions of (\ref{odePDO}) are given by
\begin{equation}
\nonumber
f_{c, \{c_k\}}=\sum_{k=0}^{\infty}a_k\varphi_{0,k}+\sum_{k,j\in \mathbb{Z}_+, j\neq j_0}
\frac{\left\langle g,\,\psi_{j,k}\right\rangle}{2^{j \alpha}+\beta} \psi_{j,k}+
\sum_{k=0}^{\infty}c_k\psi_{j_0,k}+c,
\end{equation}
where $c$, $c_k$, $k\in\z_+$, are arbitrary complex numbers, and
$a_k$, $k\in\z_+$, is a solution of system (\ref{syst}).

\end{enumerate}
\end{teo}

\textbf{Proof}
It suffices to repeat the proof of Theorem~\ref{teo_odePDO}, and take into account  that equation~(\ref{11})
with $j=j_0$ has no solution if there exists $k_0\in\z_+$ such that $\left\langle g,\,\psi_{j_0,k_0}\right\rangle\ne0$, and $\alpha_{j_0,k}=c_k$ is a solution for any $c_k\in\cn$ if   $\left\langle g,\,\psi_{j_0,k}\right\rangle=0$ for all $k\in \z_+$. \hfill $\Box$

We now consider the Cauchy problem  for the one-dimensional non-homogeneous heat equation in  variables $(x,t)$,  where $x\in G$ and  $t$ (time) is real. Let  $U=[0,\,+\infty)$  or  $U=[0,\,T].$
The problem is of the form
\begin{equation}
	\label{pdePDOheat}
	\left\{
	\begin{array}{l}
	\frac{\partial f(x,\,t)}{\partial t} = {\cal D}_x^{\alpha} f(x,\,t) + g(x,\,t), 	 \\
	f(x,\,0) = f^0(x),\\
	\end{array}
	\right.\quad x\in G,\quad t\in U.
\end{equation}

\begin{teo}
\label{teopdePDOheat}
Let $f^0 \in \tilde{S}',$ $g_t:= g(\cdot,\,t) \in \tilde{S}'$ for each $t \in U,$
and $g(x,\, \cdot)$ is   continuous on $U.$ Then
\begin{enumerate}
	\item the Cauchy problem (\ref{pdePDOheat}) has a unique solution $f(x,\,t)$ which is in $\tilde{S}'$ for every $t\in U$;
	\item this solution is in $L_2(G)$ for each $t\in U$ whenever
	$f^0  \in L_2(G)$,
	
\begin{equation}
\label{con1}
		\widehat{f^0}(\xi)=O(e^{-\|\xi\|^{\alpha} \theta({\rm log}_2 \|\xi\|)}),  \quad \|\xi\|\to \infty,
\end{equation}
	
for every $t \in U$ the distribution $g_t$ is in $L_2(G)$ and
\begin{equation}
\label{con2}
		\widehat{g_t}(\xi)=O(e^{-\|\xi\|^{\alpha} \theta({\rm log}_2 \|\xi\|)}),\quad \|\xi\|\to \infty,
\end{equation}
where  $\theta(\nu) \to \infty$ as $\nu \to \infty;$
	\item
	if either $\alpha>0$ and conditions (\ref{con1}), (\ref{con2}) are fulfilled, 	 or $\alpha<0$ and for every $t \in U$
\begin{equation}
\label{con3}
		 \widehat{f^0}(\xi)=O(\|\xi\|^{-(\varepsilon+1/2)}),\ \
\widehat{g_t}(\xi)=O(\|\xi\|^{-(\varepsilon+1/2)}),
		\quad \|\xi\|\to \infty,
\end{equation}
	where $\varepsilon>0,$
	then the  solution is continuous
on $G$ for every $t\in U$, all continuous solutions are given by $f_c(x,\,t)=f_0(x,\,t)+c,  c\in \mathbb{C}$.
\end{enumerate}
\end{teo}

\textbf{Proof.}
1.
Suppose
$f^0 = \sum_{j, k}b_{j,k}\psi_{j,k},$
$g_t = \sum_{j, k}d_{j,k}(t)\psi_{j,k}$ and
$f(\cdot,\,t) = \sum_{j, k}a_{j,k}(t)\psi_{j,k}$
are the Haar representations of $f^0$, $g_t,$ and $f(\cdot,t)$ respectively.
Using Corollary \ref{eigen_Haar}, we can rewrite~(\ref{pdePDOheat})
in the form
$$
\left\{
	\begin{array}{l}
	\sum\limits_{j\in \mathbb{Z}, k\in\mathbb{Z}_+}
\left(\dot{a}_{j,k}(t)-2^{j\alpha} a_{j,k}(t)-d_{j,k}(t)\right)\psi_{j,k}(x)=0, 	 \\
	\sum\limits_{j\in \mathbb{Z}, k\in\mathbb{Z}_+}a_{j,k}(0)\psi_{j,k}(x)=\sum\limits_{j\in \mathbb{Z}k\in\mathbb{Z}_+}b_{j,k}\psi_{j,k}(x),\\
	\end{array}
	\right.
$$
where  $\dot{a}_{j,k}$ denotes the ordinary derivative of  $a_{j,k}$ with respect to $t.$
Therefore, for every $j\in\mathbb{Z}$ and $k\in\z_+$ we obtain the  Cauchy problem for a linear ordinary differential equation of the first order
$$
\left\{
	\begin{array}{l}
	\dot{a}_{j,k}(t)-2^{j\alpha} a_{j,k}(t)-d_{j,k}(t)=0, 	\\
	a_{j,k}(0)=b_{j,k},\\
	\end{array}
	\right.
$$
 The solution to the problem is
 $$
 a_{j,k}(t)=
 \int\limits_{0}^{t} e^{2^{j \alpha}(t-\tau)}d_{j,k}(\tau)\, d\tau +
 b_{j,k}e^{2^{j \alpha}t}.
 $$
 Thus the coefficients in the Haar representation of $f$ are found.

 2.
 Let $\alpha > 0.$ The case $\alpha < 0$ is quite analogous to the first one.
 Suppose $h\in L_2(G),$ then
\begin{equation}
	\label{sum_k}
\sum_{k\in\z_+}|\left\langle h,\,\psi_{j,k}\right\rangle|^2=
\int\limits_{I_{-j-1}\setminus I_{-j}}|\widehat{h}|^2	
\end{equation}
  Indeed, let $h_j$ be the orthogonal projection of $h$ on $W_j.$ Then, taking into account   that the functions $\psi_{j,k}$ form a basis for $W_j$, using  the Parseval and Plancherel equalities, and (\ref{05}), we obtain
  $$
  \sum_{k\in\z_+}|\left\langle h,\,\psi_{j,k}\right\rangle|^2=
  \sum_{k\in\z_+}|\left\langle h_j,\,\psi_{j,k}\right\rangle|^2=
  \|h_j\|^2=\|\widehat{h_j}\|^2=\int\limits_{I_{-j-1}\setminus I_{-j}}|\widehat{h_j}|^2=\int\limits_{I_{-j-1}\setminus I_{-j}}|\widehat{h}|^2.
  $$

Fix $t\in U.$ Suppose that all conditions of item 2 are fulfilled.
Recall that
$\{\psi_{j,k}\}_{j\in \mathbb{Z},k\in\z_+}$ is an orthonormal basis in $L_2(G)$
and using the Haar representation of $f$ from item~1, we have
$$
\|f(\cdot,\,t)\|^2_{L_2(G)}=\sum_{j\in\mathbb{Z}\,k\in\z_+}|a_{j,k}|^2\leq
\sum_{j\in\mathbb{Z}\,k\in\z_+} e^{2^{j\alpha+1}t}|b_{j,k}|^2 +
\sum_{j\in\mathbb{Z}\,k\in\z_+} e^{2^{j\alpha+1}t}\left|\int\limits_{0}^{t}e^{-2^{j \alpha}\tau}d_{j,k}\, d\tau\right|^2
$$
Applying the Cauchy-Bunyakovskii inequality
and calculating the  integral
$\int_{0}^{t}e^{2^{j\alpha+1}\tau}\, d\tau$,
 we get
$$
\|f(\cdot,\,t)\|^2_{L_2(G)}\leq
\sum_{j\in\mathbb{Z}} e^{2^{j\alpha+1}t}\sum_{k\in\mathbb{Z}_+}|\left\langle f^0,\,\psi_{j,k}\right\rangle|^2
+
\sum_{j\in\mathbb{Z}}
e^{2^{j\alpha+1}t}\frac{1-e^{-2^{j\alpha+1}t}}{2^{j\alpha+1}}
\int\limits_{0}^{t}\sum_{k\in\z_+}|(g_{\tau},\,\psi_{j,k})|^2 \, d\tau
$$
Thus, by (\ref{sum_k}) with $h=f^0$ and $h=g_{\tau},$
$$
\|f(\cdot,\,t)\|^2_{L_2(G)}\leq
\sum_{j\in\mathbb{Z}} e^{2^{j\alpha+1}t}\int\limits_{I_{-j-1}\setminus I_{-j}}|\widehat{f^0}(\xi)|^2\,d\xi
+
\sum_{j\in\mathbb{Z}}
\frac{e^{2^{j\alpha+1}t}-1}{2^{j\alpha+1}}
\int\limits_{0}^{t}\int\limits_{I_{-j-1}\setminus I_{-j}}
|\widehat{g_{\tau}}(\xi)|^2\, d\xi \, d\tau.
$$

Consider 2 cases.
Case 1: $j\geq 0.$
By (\ref{con1}),
$$
\sum_{j\in\mathbb{Z}_+} e^{2^{j\alpha+1}t}\int\limits_{I_{-j-1}\setminus I_{-j}}|\widehat{f^0}(\xi)|^2\,d\xi \leq C
\sum_{j\in\mathbb{Z}_+} 2^j e^{2^{j\alpha+1}(t-\theta(j))}<\infty.
$$
By (\ref{con2}),
$$
\sum_{j\in\mathbb{Z}_+}
\frac{e^{2^{j\alpha+1}t}-1}{2^{j\alpha+1}}
\int\limits_{0}^{t}\int\limits_{I_{-j-1}\setminus I_{-j}}|\widehat{g_{\tau}}(\xi)|^2\, d\xi \, d\tau
\leq C
t \sum_{j\in\mathbb{Z}_+} \frac{(e^{2^{j\alpha+1}t}-1)e^{-2^{j\alpha+1}\theta(j)}}{2^{j(\alpha-1)+1}}.
$$
It is clear that the latter sum is finite.

Case 2: $j < 0.$
We now have
$$
\sum_{j\in\mathbb{Z}\setminus \z_+} e^{2^{j\alpha+1}t}\int\limits\limits_{I_{-j-1}\setminus I_{-j}}|\widehat{f^0}(\xi)|^2\,d\xi
+
\sum_{j\in\mathbb{Z}\setminus \z_+}
\frac{e^{2^{j\alpha+1}t}-1}{2^{j\alpha+1}}
\int\limits\limits_{0}^{t}\int\limits_{I_{-j-1}\setminus I_{-j}}|\widehat{g_{\tau}}(\xi)|^2\, d\xi \, d\tau
$$
$$
\leq
C_1(t) \sum_{j\in\mathbb{Z}\setminus \z_+} \int\limits_{I_{-j-1}\setminus I_{-j}}|\widehat{f^0}(\xi)|^2\,d\xi +
C_2(t) \sum_{j\in\mathbb{Z}\setminus \z_+} \int\limits_{0}^{t} \int\limits_{I_{-j-1}\setminus I_{-j}}|\widehat{g_{\tau}}(\xi)|^2\, d\xi \, d\tau
$$
$$
\leq
C_1(t)  \int\limits_{I}|\widehat{f^0}(\xi)|^2\,d\xi +
C_2(t) \int\limits_0 ^t \int\limits_{I}|\widehat{g_{\tau}}(\xi)|^2\, d\xi\,d\tau
\leq C_1(t) \|\widehat{f^0}\|^2 + C_2(t) \int\limits_0^t\|\widehat{g_{\tau}}\|^2\,d\tau
$$
$$
=
 C_1(t) \|f^0\|^2 + C_2(t)\int\limits_0^t\|g_{\tau}\|^2\,d\tau.
$$
It remains to note that, because of the continuity of $g(x,\,\cdot )$,
$$
\int\limits_0^t\|g_{\tau}\|^2\,d\tau=
\int\limits_G\int\limits_0^t|g({x,\tau})|^2\,d\tau\,dx=
t\int\limits_G|g({x,\tau_0})|^2\,dx, \quad 0\le\tau_0\le t,
$$
If $\alpha < 0,$ then cases 1 and 2 change places.

3.
Fix $t\in U,$ $x_0\in G$ and the compact set $I_0(x_0).$
By the definition of the Haar function,
for every $x\in I_0(x_0)$
there is a unique $k=k(x)=k(x_0)\in \mathbb{Z}_+$ such that
$x\in {\rm supp}\,\varphi_{0,k(x_0)}$, and
for every $j \in \mathbb{Z}_+$  there is a unique $k=k(j,x)\in \mathbb{Z}_+$ such that $x\in {\rm supp}\,\psi_{j,k(j,x)}.$
According to Theorem \ref{oneparam}, the exists a family of
functions $f_c(\cdot, t)\in S'$ which and only which
 have the same Haar representation as $f$.
Using the quasi-Haar representation for $f_c(\cdot, t)$ (see Theorem \ref{oneparam}), we have
$$
f_c(x,\,t)=a_{k(x)}(\,t) \varphi_{0,k(x)}(x) + \sum_{j\in \mathbb{Z}_+}a_{j,k(j,x)}(t)
\psi_{j,k(j,x)}(x)+c,
$$

Let us prove that for every $t\in U$ the series
 $ \sum_{j\in \mathbb{Z}_+}|a_{j,k(j,x)}(t)|$ converges uniformly  with respect to $x\in I_0(x_0).$
 Using the Plancherel equality and  (\ref{hat_Haar}), we get
 $$
 |b_{j,k}|=\left|\left\langle f^0,\,\psi_{j,k}\right\rangle\right|=
 2^{-j/2} \left|\int\limits_{I_{-j-1}\setminus I_{-j}}
 \widehat{f^0}(\xi) {\rm w}_{k}(D^j\xi)\,d\xi\right|\leq
 2^{-j/2} \int\limits_{I_{-j-1}\setminus I_{-j}} |\widehat{f^0}|
 $$
 and analogously
 $$
 |d_{j,k}(t)|=\left|\left\langle g_{t},\,\psi_{j,k}\right\rangle\right|\leq
 2^{-j/2} \int\limits_{I_{-j-1}\setminus I_{-j}} |\widehat{g_t}|.
 $$
 Therefore, if $\alpha > 0,$ then by (\ref{con1}) and (\ref{con2}), we have
 $$
 \sum_{j\in \mathbb{Z}_+}|a_{j,k(j,x)}(t)|\leq
\sum_{j\in\mathbb{Z}_+} e^{2^{j\alpha}t}|b_{j,k(j,x)}| +
\sum_{j\in\mathbb{Z}_+} e^{2^{j\alpha}t}\left|\int\limits_{0}^{t}e^{-2^{j \alpha}\tau}d_{j,k(j,x)}(\tau)\, d\tau\right|
$$
$$
\leq C
\sum_{j\in\mathbb{Z}_+} 2^{j/2} e^{2^{j\alpha}(t-\theta(j))}+
C
\sum_{j\in\mathbb{Z}_+} \frac{e^{2^{j\alpha}t}-1}
{2^{j (\alpha-1/2)}}e^{-2^{j\alpha}\theta(j)}.
 $$
It is clear that both the sums are finite.

If $\alpha < 0,$ then by (\ref{con3}) and elementary properties of the
exponential function,  we have
$$
 \sum_{j\in \mathbb{Z}_+}|a_{j,k(j,x)}(t)|\leq
\sum_{j\in\mathbb{Z}_+} e^{2^{j\alpha}t}|b_{j,k(j,x)}| +
\sum_{j\in\mathbb{Z}_+} e^{2^{j\alpha}t}\left|\int\limits_{0}^{t}e^{-2^{j \alpha}\tau}d_{j,k(j,x)}(\tau)\, d\tau\right|
$$
$$
\leq C
\sum_{j\in\mathbb{Z}_+} e^{2^{j\alpha}t}2^{-j\varepsilon} +
C
\sum_{j\in\mathbb{Z}_+} \frac{e^{2^{j\alpha}t}-1}{2^{j \alpha}}2^{-j\varepsilon}\leq
C e^{t}\sum_{j\in\mathbb{Z}_+}2^{-j\varepsilon} +
C \sum_{j\in\mathbb{Z}_+}\left(t+\frac{e^t t^2}{2}2^{j\alpha}\right)2^{-j\varepsilon}.
$$
It is clear that both the sums are finite.

 Therefore
$ \sum_{j\in \mathbb{Z}_+}|a_{j,k(j,x)}(t)|$ uniformly converges on $I_0(x_0).$
Thus  $f_c(\cdot, t)$ is continuous at~$x_0$.
 It remains to note that $x_0$ is an arbitrary element of $G,$
thus $f_c(\cdot, t)$ is continuous on~$G.$
\hfill $\Box$

Finally note that a statement analogous to item 1 of Theorem~\ref{teopdePDOheat}
 was proved by Shelkovich in~\cite{Sh} for the functions defined on the field of $p$-adic numbers. The idea to use the class of distributions ${\tilde S}'$ and the
Haar representation of its elements was borrowed  from~\cite{Sh}.
Several other PDEs were also considered in this paper. We
restrict ourselves to the study of the heat equation, because actually a  repetition of all arguments are needed to solve  other equations, for which this technique works.

\section*{Acknowledgments}
The work  is supported by the RFBR-grant \#12-01-00216.

\end{document}